# POINCARÉ AND TRANSPORTATION INEQUALITIES FOR GIBBS MEASURES UNDER THE DOBRUSHIN UNIQUENESS CONDITION

BY LIMING WU

*Université Blaise Pascal and Wuhan University*

In in this paper we establish an explicit and sharp estimate of the spectral gap (Poincaré inequality) and the transportation inequality for Gibbs measures, under the Dobrushin uniqueness condition. Moreover, we give a generalization of the Liggett's $M - \varepsilon$ theorem for interacting particle systems.

**1. Introduction.** Consider the configuration space $E^T$ of an interacting particle system where $E$, some Polish space, represents the spin space, and $T$, an at most countable set, represents the sites. Its equilibrium states are described by the Gibbs measures $\mu$ on $E^T$ associated with a local specification $(\mu_i = \mu_i(dx_i|x))_{i \in T}$, that is, for each $i \in T$, the conditional distribution $\mu(\cdot/x)$ of $x_i$ knowing $x_{T\setminus\{i\}}$ coincides with the given $\mu_i(\cdot|x)$.

In the free (no interaction) case, $\mu_i$ is independent of $x_{T\setminus\{i\}}$ and $\mu = \prod_{i \in T} \mu_i$. In that case, we have the following Efron–Stein inequality:

$$(1.1) \qquad \lambda_1(\mu)\mu(f,f) \leq \mathbb{E}^\mu \sum_{i \in T} \mu_i(f,f) \qquad \forall f \in L^2(E^T, \mu),$$

where $\lambda_1(\mu) = 1$, $\mu(f,g)$ denotes the covariance of $f,g$ under $\mu$, and $\mu_i(f,g) = \mu_i(fg) - \mu_i(f)\mu_i(g)$ is the conditional covariance of $f,g$ under $\mu_i$ with $x_{T\setminus\{i\}}$ fixed. That $\lambda_1(\mu) = 1$ is sharp can be seen for functions $f(x) = f(x_i)$. This inequality is very important for concentration inequalities in statistics and statistical learning; see the St Flour course by Massart [26].

Our objective is to generalize this inequality to a Gibbs measure with interaction. A first crucial idea is to interpret (1.1) as a Poincaré inequality.

To that end, consider the generator

(1.2) $$\mathcal{L}f := \sum_{i \in T}[\mu_i(f) - f].$$

It generates a Glauber dynamics (a well-known stochastic algorithm) which is a Markov process of pure jumps described intuitively as follows: if the configuration at present is $x$, then at each site $i$, it will change to $y_i$ according to the distribution $\mu_i(dy_i|x)$ at rate 1 (which depends only on $x_{T \setminus \{i\}}$, not on $x_i$). The semigroup $(P_t)$ generated by $\mathcal{L}$ is symmetric on $L^2(\mu)$ and the associated Dirichlet form is exactly

$$\mathcal{E}(f, f) = \mathbb{E}^\mu \sum_{i \in T} \mu_i(f, f).$$

So inequality (1.1) becomes a Poincaré inequality for $\mathcal{L}$, and it means that $P_t$ decays to the equilibrium measure $\mu$ exponentially at speed $e^{-\lambda_1(\mu)t}$, in $L^2(\mu)$.

When $E$ is a continuous spin space, that is, a connected and complete Riemannian manifold, it is more natural to estimate the spectral gap $\lambda_1(\mu, \nabla)$, that is, the best constant in the following Poincaré inequality:

(1.3) $$\lambda_1(\mu, \nabla)\mu(f, f) \leq \mathbb{E}^\mu \sum_{i \in T} |\nabla_i f|^2,$$

where $\nabla_i$ is the gradient acting on the $i$th variable $x_i \in E$.

The studies on $\lambda_1(\mu)$ and $\lambda_1(\mu, \nabla)$ (and related inequalities) are part of a long story and the field remains very active. Let us review a series of works which motivate directly our study (the reader is referred to [16, 23] for numerous related references).

The first important and general result of quantitative type is Liggett's $M - \varepsilon$ theorem, which gives an explicit exponential decay $e^{-(\varepsilon - M)t}$ of a general interacting particle system (of pure jumps) to its equilibrium measure, in the *triple norm of Liggett* (see [21], Chapter I, Theorem 3.8 for the definition of $\varepsilon$ and $M$). Applied to $\mathcal{L}$ given by (1.2), it yields $\lambda_1(\mu) \geq \varepsilon - M$ by Lemma 2.6 in this paper. When $E$ has exactly two elements, Liggett's $M - \varepsilon$ theorem for $\mathcal{L}$ may be regarded as a dynamical counterpart of Dobrushin's uniqueness criterion [8, 9]. But Liggett's estimate is no longer accurate when $E$ has more than two elements, and becomes inapplicable for *infinite* spin space $E$ (since $\varepsilon = 0$ in such case). Of course Liggett's theorem does not furnish information about $\lambda_1(\mu, \nabla)$.

Recall that in the two-points spin space case, Maes and Shlosman [22] have found a constructive criterion for the validity of Liggett's exponential convergence (but without an explicit estimate better than Liggett's theorem), which becomes necessary for the attractive system.

In an important contribution [38], Zegarlinski proved the logarithmic Sobolev (log-Sobolev) inequality, which is stronger than the Poincaré inequality, with an explicit constant under some condition which is inspired by the Dobrushin uniqueness condition, both for continuous spin space or two-points spin space. His condition (see (0.12) and (0.13) in [38]), though in spirit quite close to the Dobrushin uniqueness condition (see Section 2), is in reality very different, as seen for a number of concrete examples (already discussed in [38]; see also Section 5).

The most spectacular advances were made on the qualitative aspect of $\lambda_1(\mu)$ and $\lambda_1(\mu, \nabla)$, that is, about the validity of the Poincaré inequality and the log-Sobolev inequality. When the spin space $E$ is finite or compact, $T$ is the lattice $\mathbb{Z}^d$, and $(\mu_i)$ is given by a family of interaction functions $(\phi_S)$ with finite range, Stroock and Zegarlinski [30, 31] prove essentially the equivalence between the Poincaré inequality, the log-Sobolev inequality and the Dobrushin–Shlosman complete analyticity (CA) condition (see [10, 30, 31] for CA). Their method for establishing those inequalities consists in an iteration procedure, which does not provide explicit estimates of the involved constants, unlike [38]. See [23] and [16] for further development and references. The (partial) extension of their impressive results to the unbounded spin case for Glauber dynamics with a single-site diffusion term was carried out by Helffer [17], Ledoux [19] and Yoshida [37] and so on. In particular Helffer [17] gave some explicit estimates of $\lambda_1(\mu, \nabla)$ by means of the Witten Laplacian, and Ledoux [19] realized it by a very simple and elegant argument based on the $\Gamma_2$-technique. More recently, the author [36] obtained an explicit and sharp estimate for a continuous gas, based on Liggett's $M - \varepsilon$ theorem.

The advantage of the Dobrushin uniqueness condition over the Dobrushin–Shlosman CA is the following: (1) the Dobrushin uniqueness condition is quantitative and explicit; (2) it holds for general graphs $T$ and general interaction of infinite range. (But in the lattice and finite range case, the Dobrushin uniqueness condition is more restrictive than the Dobrushin–Shlosman CA.)

In this paper we shall not only provide some explicit sharp estimates of both $\lambda_1(\mu)$ and $\lambda_1(\mu, \nabla)$ based directly on the Dobrushin uniqueness condition, but also present a unified and particularly simple approach by generalizing Liggett's $M - \varepsilon$ theorem (avoiding so the high technical difficulties existing in the known works mentioned above).

Indeed, we shall derive the estimate of $\lambda_1(\mu, \nabla)$ from that of $\lambda_1(\mu)$. Our approach for estimating $\lambda_1(\mu)$ is based on an exchange relation between $\mathcal{L}$ and the difference operator $D_{x_j \to y_j}$, largely inspired by the Bochner formula for the commutator between the Laplacian and the gradient on Riemannian manifold or the $\Gamma_2$-method of Bakry–Emery–Ledoux. In reality this method allows us to generalize Liggett's $M - \varepsilon$ theorem from the single-site *finite* spin space to a general (possibly continuous or unbounded) single-site spin space, which is stronger than the Poincaré inequality.

This paper is organized as follows. In Section 2 we recall Dobrushin's interdependence coefficients of a Gibbs measure and establish a sharp lower bound for $\lambda_1(\mu)$ and $\lambda_1(\mu, \nabla)$ under the Dobrushin uniqueness condition. Our method described above can be easily generalized to general interacting particle systems, leading to an extension of Liggett's $M - \varepsilon$ theorem. This is carried out in Section 3.

In Section 4 we investigate another object of this paper: the $L^1$-transportation inequality for the Gibbs measure. Indeed, we shall interpret the famous Dobrushin a priori estimate as a variant of (1.1) for the Wasserstein metric. From that new version we derive easily the $L^1$-transportation inequality by the martingale method and get Hoeffding's Gaussian concentration inequality as a corollary. This extends the corresponding work of Marton [24] and Djellout, Guillin and Wu [7] on contracting Markov chains.

Several concrete examples are provided in Section 5 for illustrating our general results.

**2. Poincaré inequality for Gibbs measure.** As we are mainly interested in the explicit estimate on the spectral gap for Gibbs measures, so it is enough to get such an estimate in finite volume, independent of boundary condition and of the finite volume. That is why we shall work with a probability measure $\mu$ on $E^T$ only with $T$ finite, throughout this paper.

Throughout this paper $(E, \mathcal{B}, d)$ is fixed as follows: either $(E, d)$ is a metrical complete and separable (say, Polish) space equipped with the Borel field $\mathcal{B}$; or $(E, \mathcal{B})$ is a measurable space and $d(x, y) = \mathbb{1}_{x \neq y}$ (the trivial metric) such that $d$ is $\mathcal{B} \times \mathcal{B}$-measurable.

2.1. *Dobrushin's interdependence coefficients.* Let $M_1(E)$ be the space of probability measures on $(E, \mathcal{B})$ and $M_1^{d,p}(E) := \{\nu \in M_1(E); (\int_E d^p(x_0, x) \times \nu(dx))^{1/p} < +\infty\}$ ($x_0 \in E$ is some fixed point), where $1 \leq p \leq +\infty$. Given $\nu_1, \nu_2 \in M_1^{d,p}(E)$, the $L^p$-Wasserstein distance between $\nu_1, \nu_2$ is given by

$$(2.1) \qquad W_{p,d}(\nu_1, \nu_2) := \inf_\pi \left( \iint_{E \times E} d(x, y)^p \pi(dx, dy) \right)^{1/p},$$

where the infimum is taken over all probability measures $\pi$ on $E \times E$ such that its marginal distributions are respectively $\nu_1$ and $\nu_2$. When $d(x, y) = \mathbb{1}_{x \neq y}$ (the trivial metric), it is well known that

$$W_{1,d}(\nu_1, \nu_2) = \sup_{A \in \mathcal{B}} |\nu_1(A) - \nu_2(A)| = \tfrac{1}{2} \|\nu_1 - \nu_2\|_{\mathrm{TV}} \qquad \text{(total variation)}.$$

Recall (cf. [32])

$$(2.2) \quad W_{1,d}(\mu, \nu) = \sup_{\|f\|_{\mathrm{Lip}} \leq 1} \int_E f \, d(\mu - \nu), \qquad \|f\|_{\mathrm{Lip}} := \sup_{x \neq y} \frac{|f(x) - f(y)|}{d(x, y)}.$$

Let $\mu_i(dx_i|x)$ be the given regular conditional distribution of $x_i$ knowing $x_{T\setminus\{i\}}$. Define the $d$-Dobrushin interdependence matrix $C := (c_{ij})_{i,j \in T}$ by

$$(2.3) \qquad c_{ij} := \sup_{x=y \text{ off } j} \frac{W_{1,d}(\mu_i(\cdot/x), \mu_i(\cdot/y))}{d(x_j, y_j)} \qquad \forall i, j \in T$$

(obviously $c_{ii} = 0$). Then the Dobrushin uniqueness condition [8, 9] is

$$\sup_i \sum_j c_{ij} < 1.$$

Notice that the l.h.s. above coincides with the norm $\|C\|_\infty$ of $C : l_\infty(T) \to l_\infty(T)$.

2.2. *Sharp estimates of the spectral gap.* The main result of this section is the following:

THEOREM 2.1. *Assume that $\int_{E^T} d(x_i, y_i)^2 \, d\mu(x) < +\infty$ for all $i \in T$ and for some (fixed) $y \in E^T$. Let $r_{\mathrm{sp}}(C)$ be the spectral radius of the Dobrushin matrix $C = (c_{ij})_{i,j \in T}$ (which is an eigenvalue of $C$ by the Perron–Frobenius theorem).*

*If $r_{\mathrm{sp}}(C) < 1$, then*

$$(2.4) \qquad (1 - r_{\mathrm{sp}}(C))\mu(f, f) \leq \mathbb{E}^\mu \sum_{i \in T} \mu_i(f, f) \qquad \forall f \in L^2(E^T, \mu).$$

*In particular, the lowest eigenvalue above zero $\lambda_1(\mu)$ of $-\mathcal{L}$ in $L^2(E^T, \mu)$ (called the spectral gap of $\mu$) verifies*

$$\lambda_1(\mu) \geq 1 - r_{\mathrm{sp}}(C),$$

*where $\mathcal{L}f := \sum_{i \in T}[\mu_i(f) - f] \; \forall f \in L^2(E^T, \mu)$.*

It is an elementary fact that $r_{\mathrm{sp}}(C) \leq \|C\|_\infty = \sup_i \sum_j c_{ij}$, the quantity in the Dobrushin uniqueness condition. In the free case [i.e., $\mu(\cdot|x)$ is independent of $x$] $C = 0$ and the inequality (2.4) becomes the sharp Efron–Stein inequality (1.1). Notice that the estimate (2.4) on the spectral gap above depends sensitively on the choice of the metric $d$ via the Dobrushin interdependence matrix $C$, which allows us to apply it for the discrete or continuous, compact or noncompact spin spaces.

When $E$ is a complete and connected Riemannian manifold equipped with the Riemannian metric $d$, the following pre-Dirichlet form,

$$(2.5) \qquad \mathcal{E}^\nabla(f, f) := \mathbb{E}^\mu \sum_{i \in T} |\nabla_i f|^2 \qquad \forall f \in C_b^1(E^T),$$

is more often used to generate the Glauber dynamics, where $\nabla_i$ is the gradient acting on the $i$th coordinate $x_i \in E$. Let $\lambda_0$ be the infimum of the spectral gap of $\mu_i(\cdot/x)$ w.r.t. $\nabla_i$ on $E$ over all $i \in T$ and $x \in E^T$, more precisely,

(2.6)
$$\lambda_0 := \sup\left\{\lambda \geq 0 | \lambda \mu_i(f,f) \leq \int_E |\nabla_i f|^2 \, d\mu_i \right.$$
$$\left. \forall x \in E^T, i \in T, \ \forall f \in C_b^1(E^T)\right\}.$$

Assume that $\lambda_0 > 0$, then with $\mathcal{E}^\nabla(f,f)$ given by (2.5),

$$\mathbb{E}^\mu \sum_{i \in T} \mu_i(f,f) \leq \frac{1}{\lambda_0} \mathcal{E}^\nabla(f,f) \qquad \forall f \in C_b^1(E^T).$$

So we derive immediately from Theorem 2.1 the following:

THEOREM 2.2.  *In the context above, assume that $\int_{E^T} d(x_i, y_i)^2 \, d\mu(x) < +\infty$ for all $i \in T$ and for some (fixed) $y \in E^T$. If $r_{\mathrm{sp}}(C) < 1$ and $\lambda_0 > 0$, then*

(2.7) $\qquad \lambda_0(1 - r_{\mathrm{sp}}(C))\mu(f,f) \leq \mathbb{E}^\mu \sum_{i \in T} |\nabla_i f|^2 \qquad \forall f \in C_b^1(E^T),$

*that is, the best constant $\lambda_1(\mu, \nabla)$ for the Poincaré inequality (1.3) verifies*

$$\lambda_1(\mu, \nabla) \geq \lambda_0(1 - r_{\mathrm{sp}}(C)).$$

Obviously in the free case this inequality is sharp, just as the inequality (2.4). The reader might wonder if this procedure of deriving (2.7) from Theorem 2.1 is sharp in the dependent case. The following example shows that both Theorems 2.1 and 2.2 are sharp in the dependent case.

EXAMPLE 2.3 (Gaussian model). Let $T = \{1,2\}$, $E = \mathbb{R}$ and $\mu$ be the centered Gaussian measure on $\mathbb{R}^2$ such that

$$\mu(x_i, x_i) = 1, \qquad i = 1,2; \qquad \mu(x_1, x_2) = \rho \neq 0.$$

We note the following:

(i) If $\{i,j\} = \{1,2\}$, the conditional law $\mu(dx_i|x)$ of $x_i$ knowing $x_j$ is the Gaussian law $\mathcal{N}(m, \sigma^2)$ with

$$m = \rho x_j, \qquad \sigma^2 = 1 - \rho^2.$$

(ii) Since $W_{p,d}(N(m_1, \sigma^2), N(m_2, \sigma^2)) = |m_1 - m_2|$ for all $1 \leq p \leq +\infty$ (left to the reader), then the Dobrushin coefficients w.r.t. the Euclidean metric are given by

$$c_{12} = c_{21} = |\rho|.$$

Hence, $r_{\mathrm{sp}}(C) = |\rho|$.

- *Sharpness of Theorem* 2.1. By (2.4), $\lambda_1(\mu) \geq 1 - |\rho|$. We claim that $\lambda_1(\mu) = 1 - |\rho|$, showing the sharpness of Theorem 2.1. In fact, taking $f(x_1, x_2) := x_1 + (\rho/|\rho|)x_2$, we have

$$\lambda_1(\mu) \leq \frac{\mathbb{E}^\mu[\mu_1(f,f) + \mu_2(f,f)]}{\mu(f,f)} = \frac{2(1-\rho^2)}{2(1+|\rho|)} = 1 - |\rho|.$$

- *Sharpness of Theorem* 2.2. Recall at first that the Gaussian measure $\nu = \mathcal{N}(m, \Gamma)$ on $\mathbb{R}^d$ with mean $m \in \mathbb{R}^d$ and covariance matrix $\Gamma$ satisfies

$$\lambda_1(\nu, \nabla) = \frac{1}{\lambda_{\max}(\Gamma)},$$

where $\lambda_{\max}(\Gamma)$ is the maximal eigenvalue of $\Gamma$.

The covariance matrix of $\mu$ is given by $\Gamma = \begin{pmatrix} 1 & \rho \\ \rho & 1 \end{pmatrix}$. Its maximal eigenvalue is $1 + |\rho|$. Hence, $\lambda_1(\mu, \nabla) = (1 + |\rho|)^{-1}$.

Let us see why Theorem 2.2 produces the same result. In fact, the spectral gap of $\nu := N(m, \sigma^2)$ w.r.t. the Dirichlet form $\int_\mathbb{R} f'^2 \, d\nu$ is exactly $1/\sigma^2$. Hence, $\lambda_0$ defined in (2.6) equals $(1 - \rho^2)^{-1}$. Hence, the lower bound of $\lambda_1(\mu, \nabla)$ given in Theorem 2.2 becomes

$$\lambda_0(1 - r_{\mathrm{sp}}(C)) = \frac{1 - |\rho|}{1 - \rho^2} = \frac{1}{1 + |\rho|},$$

which is the exact value of $\lambda_1(\mu, \nabla)$. Thus, for this example, Theorem 2.2 is sharp.

REMARK 2.4. Let us compare our results with the explicit estimate of the constant in the log-Sobolev inequality by Zegarlinski [38]. In the continuous spin space case, his assumptions are the following:

(Z1) There is some $c_0 \in (0, \infty)$ such that

$$Ent_{\mu_i}(f^2) \leq 2c_0 \int_E |\nabla_i f|^2 \, d\mu_i \qquad \forall i \in T, x \in E^T, f \in C_b^1(E).$$

Here $Ent_\nu(f) := \nu(f \log f) - \nu(f) \log \nu(f)$ is the Kullback entropy of $0 \leq f \in L^1(\nu)$ w.r.t. the probability measure $\nu$.

(Z2) There exist $C^Z = (c_{ij}^Z \geq 0)_{i,j \in T}$ such that

(2.8) $$|\nabla_j(\mu_i(f^2))^{1/2}| \leq [\mu_i(|\nabla_j f|^2)]^{1/2} + c_{ij}^Z[\mu_i(|\nabla_i f|^2)]^{1/2}$$
$$\forall f \in C_b^1(E^T)$$

and

$$\gamma = \max(\|C^Z\|_1, \|C^Z\|_\infty) < 1.$$

Under those conditions, he derived the following log-Sobolev inequality:

$$Ent_{\mu_T}(f^2) \leq \frac{2c_0}{(1-\gamma)^2}\mathcal{E}^\nabla(f,f) \qquad \forall f \in C_b^1(E^T). \tag{2.9}$$

It implies that $\lambda_1(\mu_T, \nabla) \geq (1-\gamma)^2/c_0$. But this estimate is in general less accurate than Theorem 2.2, for $c_{ij}^Z$ is much more difficult to compute and is in general much larger than $c_{ij}$ for concrete models (see discussions in Section 5). Another important advantage of our approach is that we can choose a metric (not necessarily the Riemannian) w.r.t. which our $c_{ij}$ becomes as small as possible. This will be illustrated in Example 5.3.

In the two-points spin space $E = \{+, -\}$ case, the explicit estimate of the constant in the log-Sobolev inequality by Zegarlinski [38], Theorem 4.3, is much more larger and works only in the finite range case.

2.3. *Proof of Theorem* 2.1. We shall prove it from a (Glauber) dynamical point of view, that is, via the analysis of the semigroup $P_t = e^{t\mathcal{L}}$. Let $\mathcal{F}_T := \mathcal{B}^T$ and $b\mathcal{F}_T$ the Banach space of real bounded and measurable functions on $(E^T, \mathcal{F}_T)$ equipped with the sup norm. As $\mathcal{L}f = \sum_{i \in T}[\mu_i(f) - f]$ is bounded on $b\mathcal{F}_T$, $P_t = e^{t\mathcal{L}}$ is a well-defined Markov semigroup on $b\mathcal{F}_T$.

Consider the space of Lipschitzian continuous functions $C_{\text{Lip}}(E^T) := \{f \in b\mathcal{F}_T; \sup_{i \in T} \delta_i(f) < +\infty\}$, where

$$\delta_i(f) = \sup_{x=y \text{ off } i} \frac{|f(y) - f(x)|}{d(y_i, x_i)}.$$

In fact, we shall prove the stronger.

PROPOSITION 2.5. *In the context of Theorem* 2.1, *let* $P_t = e^{t\mathcal{L}}$. *Then for all measurable* $f \in b\mathcal{F}_T$ *such that* $\sup_{i \in T} \delta_i(f) < +\infty$,

$$\delta_j(P_t f) \leq \sum_{i \in T} \delta_i(f)(e^{-t(I-C)})_{ij} \qquad \forall j \in T, \tag{2.10}$$

*where* $I$ *is the identity matrix. In particular,*

$$\sum_{j \in T} \delta_j(P_t f) \leq e^{-t(1-\|C\|_\infty)} \sum_{j \in T} \delta_j(f),$$

$$\max_{j \in T} \delta_j(P_t f) \leq e^{-t(1-\|C\|_1)} \max_{j \in T} \delta_j(f), \tag{2.11}$$

*where*

$$\|C\|_\infty := \sup_i \sum_j c_{ij}; \qquad \|C\|_1 := \|C\|_{l_1(T) \to l_1(T)} = \sup_j \sum_i c_{ij}.$$

Once this result is proved, Theorem 2.1 follows immediately from (2.10) and the following general fact [by choosing $\mathcal{D} = C_{\text{Lip}}(E^T)$]:

LEMMA 2.6 ([34, 35]). *Let $(P_t)$ be a strongly continuous symmetric Markov semigroup on $L^2(\mu)$. Assume that there are a dense subset $\mathcal{D} \subset L^2(\mu)$ and a constant $\delta > 0$ such that $\forall f \in \mathcal{D}$, $\exists C(f) > 0$,*

$$\mu(P_t f, P_t f) \leq C(f) e^{-2\delta t} \qquad \forall t > 0,$$

*then the spectral gap $\lambda_1$ of $-\mathcal{L}$ verifies $\lambda_1 \geq \delta$.*

We now turn to the proof.

PROOF OF PROPOSITION 2.5. Notice at first that, for the transposition $C^t$, we have

$$\|C^t\|_1 = \|C\|_\infty, \qquad \|C^t\|_\infty = \|C\|_1.$$

Hence, (2.11) follows from (2.10). Below we prove (2.10) in three steps.

*Step* 1. Consider the difference operator $D_{x_j \to y_j} f(x) := f(x^{y_j}) - f(x)$, where

$$x^{y_j}(i) := \begin{cases} x(i), & \text{if } i \neq j, \\ y_j, & \text{if } i = j, \end{cases}$$

and $y_j$ is some fixed point of $E$. Our idea is based, roughly, say, on the calculation of $D_{x_j \to y_j} \mathcal{L} f - \mathcal{L} D_{x_j \to y_j} f$, where $f \in C_{\text{Lip}}(E^T)$. For such $f$, when $i = j$, we have

$$D_{x_j \to y_j}[\mu_j(f) - f] = -D_{x_j \to y_j} f,$$

for $\mu_j(f)$ is $\mathcal{F}_{T\setminus\{j\}}$-measurable. When $i \neq j$, by putting $x^{z_i y_j}(k) = z_i$ if $k = i$ and $y_j$ if $k = j$ and $x_k$ otherwise, we have

$$D_{x_j \to y_j}[\mu_i(f) - f](x)$$
$$= \int_E \mu_i(dz_i|x^{y_j})(f(x^{z_i y_j}) - f(x^{y_j})) - \int_E \mu_i(dz_i|x)(f(x^{z_i}) - f(x))$$
$$= \int_E (\mu_i(dz_i|x^{y_j}) - \mu_i(dz_i|x))(f(x^{z_i y_j}) - f(x^{y_j}))$$
$$+ \int_E \mu_i(dz_i|x)((D_{x_j \to y_j} f)(x^{z_i}) - (D_{x_j \to y_j} f)(x)).$$

Since the Lipschitzian coefficient of $z_i \to f(x^{z_i y_j}) - f(x^{y_j})$ is not greater than $\delta_i(f)$, by (2.2) and the definition of $c_{ij}$, the function

$$g_{ij}(x) := \int_E (\mu_i(dz_i|x^{y_j}) - \mu_i(dz_i|x))(f(x^{z_i y_j}) - f(x^{y_j}))$$

verifies

$$|g_{ij}| \leq c_{ij} \delta_i(f) d(x_j, y_j).$$

Thus, we obtain the following relation about $[D_{x_j \to y_j}, \mathcal{L}]$:

$$D_{x_j \to y_j} \mathcal{L} f = -D_{x_j \to y_j} f$$
(2.12)
$$+ \sum_{i:\, i \neq j} [\mu_i(D_{x_j \to y_j} f) - D_{x_j \to y_j} f] + \sum_{i:\, i \neq j} g_{ij},$$

$$|g_{ij}| \leq c_{ij} \delta_i(f) d(x_j, y_j).$$

Let $\mathcal{L}_j f := \sum_{i:\, i \neq j} [\mu_i(f) - f]$. Dividing both sides of (2.12) by $d(x_j, y_j)$ (with the convention that $0/0 := 0$) and noting the following obvious but crucial relation

$$\frac{\mathcal{L}_j g(x)}{d(x_j, y_j)} = \left( \mathcal{L}_j \frac{g}{d(\cdot_j, y_j)} \right)(x),$$
(2.13)

we get

$$\delta_j(\mathcal{L} f) \leq (\|\mathcal{L}_j\|_{b\mathcal{F}_T} + 1)\delta_j(f) + \sum_{i:\, i \neq j} c_{ij} \delta_i(f).$$

Hence, $f \to \mathcal{L} f$ is bounded on $C_{\mathrm{Lip}}$.

*Step* 2. Note that $\mathcal{L}_j f := \sum_{i:\, i \neq j} [\mu_i(f) - f]$ generates again a Markov semigroup on $b\mathcal{F}_T$. For any $\lambda > 0$ and $f \in C_{\mathrm{Lip}}(E^T)$, let $g := f - \lambda \mathcal{L} f$. By (2.12), we have

$$(\lambda + 1) D_{x_j \to y_j} f - \lambda \mathcal{L}_j D_{x_j \to y_j} f = D_{x_j \to y_j} g + \lambda \sum_{i:i \neq j} g_{ij}.$$

Dividing both sides by $d(x_j, y_j)$ and putting $h_j(x) := \frac{D_{x_j \to y_j} f(x)}{d(x_j, y_j)}$, we obtain, by (2.12) and (2.13),

$$(\lambda + 1) h_j - \lambda \mathcal{L}_j h_j \leq \delta_j(g) + \lambda \sum_{i:\, i \neq j} c_{ij} \delta_i(f).$$

As the resolvent $(\lambda + 1 - \lambda \mathcal{L}_j)^{-1}$ on $b\mathcal{F}_T$ is positive with norm bounded by $(\lambda + 1)^{-1}$, we get

$$\sup_x h_j(x) \leq \frac{1}{\lambda + 1} \left( \delta_j(g) + \lambda \sum_{i:\, i \neq j} c_{ij} \delta_i(f) \right).$$

Since $y_j$ is arbitrary, we get

(2.14)
$$\delta_j(f) \leq \frac{1}{\lambda + 1} \left( \delta_j(g) + \lambda \sum_{i:\, i \neq j} c_{ij} \delta_i(f) \right).$$

*Step* 3. With the crucial estimate (2.14) in hand, the rest of the proof is routine by following Liggett [21], Chapter I, Theorem 3.8. Now $\mathcal{L}$ is bounded on $C_{\mathrm{Lip}}(E^T)$, hence, for all $\lambda > 0$ sufficiently small, say, $\lambda \in (0, \lambda_0)$, where $0 < \lambda_0 < 1/r_{\mathrm{sp}}(C)$, $(1 - \lambda \mathcal{L})^{-1}$ is bounded on $C_{\mathrm{Lip}}(E^T)$ and $(1 - \lambda C)^{-1}$ is bounded on $\mathbb{R}^T$. On $\mathbb{R}^T$ considering the partial order $u \leq v$ iff $u_i \leq v_i$ for all $i$, and regarding $\delta(f) := (\delta_i(f))_{i \in T}$ as column vector in $\mathbb{R}^T$, we have, by (2.14),

$$\left(1 - \frac{\lambda}{\lambda+1} C^t\right) \delta((1 - \lambda \mathcal{L})^{-1} g) \leq \frac{1}{\lambda+1} \delta(g) \qquad \forall g \in C_{\mathrm{Lip}}(E^T).$$

Since $(1 - \lambda C)^{-1} = \sum_{n=0}^{\infty} \lambda^n C^n$ is a positive matrix, we get

$$\delta((1 - \lambda \mathcal{L})^{-1} g) \leq (\lambda + 1 - \lambda C^t)^{-1} \delta(g)$$

and, consequently, for all $n \in \mathbb{N}^*$,

$$\delta((1 - \lambda \mathcal{L})^{-n} g) \leq (\lambda + 1 - \lambda C^t)^{-n} \delta(g).$$

Finally, for each $t > 0$, letting $\lambda = t/n$ and $n$ go to infinity in the above estimate, we obtain

$$\delta(P_t g) \leq e^{-t(I - C^t)} \delta(g),$$

the desired estimate (2.10). □

**3. Generalization of Liggett's $M - \varepsilon$ theorem.** In this section we consider the more general generator on $E^T$ given by

$$(3.1) \qquad \mathcal{L}f(x) := \sum_{S \subset T} \int_{E^S} J_S(x, dz_S)[f(x^{z_S}) - f(x)] \qquad \forall f \in C_{\mathrm{Lip}}(E^T),$$

where the (local) jump rate $J_S(x, dz_S)$ is a bounded nonnegative kernel on $E^T \times \mathcal{F}_S$. Assume that for each $S \subset T$ and for every $j \in T$, there is some finite optimal constant $c_S(j) \geq 0$,

$$(3.2) \quad \begin{aligned} \sup_{x=y \text{ off } j} \frac{1}{d(x_j, y_j)} \bigg| \int_{E^S} g(z_S)(J_S(x, dz_S)) - J_S(y, dz_S) \bigg| \\ \leq c_S(j) \sum_{i \in S} \delta_i(g) \end{aligned}$$

for all $g \in C_{\mathrm{Lip}}(E^T)$. Note that if $d$ is the trivial metric,

$$c_S(j) \leq \tfrac{1}{2} \sup_{x=y \text{ off } j} \|J_S(x, \cdot) - J_S(y, \cdot)\|_{\mathrm{TV}}.$$

The following generalizes Proposition 2.5.

THEOREM 3.1. *Let*

$$(3.3) \qquad c_{ij} := \sum_{S \ni i} c_S(j)$$

*and*

$$\eta := \inf_{x \in E^T} \inf_{i \in T} \sum_{S \ni i} J_S(x, E^S). \tag{3.4}$$

Then $P_t = e^{t\mathcal{L}}$ is a Markov semigroup on $b\mathcal{F}_T$, mapping $C_{\text{Lip}}(E^T)$ into itself, such that, for any $f \in C_{\text{Lip}}(E^T)$,

$$\delta_j(P_t f) \leq \sum_{i \in T} \delta_i(f)(e^{-t(\eta - C)})_{ij} \qquad \forall j \in T. \tag{3.5}$$

REMARK 3.2. When $J_S(x, dz_S) = \mu_i(dz_i|x)$ for $S = \{i\}$ and 0 otherwise, this result becomes exactly Proposition 2.5. In the symmetric case, by Lemma 2.6, the estimate (3.5) in Theorem 3.1 implies that the spectral gap $\lambda_1(\mathcal{L})$ of $\mathcal{L}$ in $L^2(\mu)$ satisfies

$$\lambda_1(\mathcal{L}) \geq \eta - r_{\text{sp}}(C).$$

REMARK 3.3. There is a quite subtle point in this result: if $(\tilde{J}_S)$ is another family of jump rates such that $\tilde{J}_S(x, \cdot \setminus \{x_S\}) = J_S(x, \cdot \setminus \{x_S\})$, then $(\tilde{J}_S)$ determines the same generator $\mathcal{L}$ as $(J_S)$. One must choose $J_S(x, \{x_S\})$ so that $J_S(x, E^S)$ is independent of $x$ for ensuring $C_S(j) < +\infty$. For instance, in the framework of Proposition 2.5, $J_S = \mu_i$ for $S = \{i\}$ seems to be the best choice.

REMARK 3.4. Let us compare this result with Liggett's $M - \varepsilon$-theorem (cf. [21], Chapter I, Theorem 3.8) and the results of Maes–Shlosman [22]:

(i) At first for continuous spins, the constant $\varepsilon$ in Liggett's $M - \varepsilon$ theorem becomes zero so that it cannot be applied for obtaining the exponential convergence in that situation.

(ii) In Liggett's $M - \varepsilon$ theorem, only the trivial metric is used and its proof seems not to work for more general metrics.

(iii) The method of Maes–Shlosman [22], based on the coupling method and a time-discretization procedure, works well for two-states spin space $E$, but seems difficult to work in the present general setting.

PROOF OF THEOREM 3.1. The proof is similar to that of Proposition 2.5, but with an important difference in Step 2.

*Step* 1. Fix $f \in C_{\text{Lip}}(E^T)$, $j \in T$ and $y_j \in E$. When $j \notin S$, by putting $x^{z_S y_j}(k) = z_i$ if $k \in S$, $y_j$ if $k = j$ and $x_k$ otherwise, we have

$$D_{x_j \to y_j} \left[ \int_{E^S} J_S(\cdot, dz_S) f(\cdot^{z_S}) - f(\cdot) \right](x)$$
$$= \int_{E^S} J_S(x^{y_j}, dz_S)(f(x^{z_S y_j}) - f(x^{y_j}))$$

$$-\int_{E^S} J_S(x, dz_S)(f(x^{z_S}) - f(x))$$

$$= \int_{E^S} (J_S(x^{y_j}, dz_S) - J_S(x, dz_S))(f(x^{z_S y_j}) - f(x^{y_j}))$$

$$+ \int_{E^S} J_S(x, dz_S)((D_{x_j \to y_j} f)(x^{z_S}) - (D_{x_j \to y_j} f)(x)).$$

Since $\delta_i[z_i \to f(x^{z_S y_j}) - f(x^{y_j})] \leq \delta_i(f)$ for each $i \in S$, the function

$$g_{Sj}(x) := \int_{E^S} (J_S(x^{y_j}, dz_S) - J_S(x, dz_S))(f(x^{z_S y_j}) - f(x^{y_j}))$$

verifies, by the definition of $c_S(j)$ in (3.2),

$$|g_{Sj}| \leq c_S(j) d(x_j, y_j) \sum_{i \in S} \delta_i(f).$$

Now letting $j \in S$, we have

$$D_{x_j \to y_j} \left[ \int_{E^S} J_S(\cdot, dz_S) f(\cdot^{z_S}) - f(\cdot) \right](x)$$

$$= \int_{E^S} J_S(x^{y_j}, dz_S)(f(x^{z_S}) - f(x^{y_j})) - \int_{E^S} J_S(x, dz_S)(f(x^{z_S}) - f(x))$$

$$= \int_{E^S} (J_S(x^{y_j}, dz_S) - J_S(x, dz_S))(f(x^{z_S}) - f(x^{y_j}))$$

$$- \int_{E^S} J_S(x, dz_S)((D_{x_j \to y_j} f)(x))$$

$$= g_{Sj}(x) - J_S(x, E^S) D_{x_j \to y_j} f(x),$$

where

$$g_{Sj}(x) := \int_{E} (J_S(x^{y_j}, dz_S) - J_S(x, dz_S))(f(x^{z_S}) - f(x^{y_j}))$$

satisfies again [by the definition of $c_S(j)$ in (3.2)]

$$|g_{Sj}| \leq c_S(j) d(x_j, y_j) \sum_{i \in S} \delta_i(f).$$

Thus, letting

$$\mathcal{L}_j f(x) := \sum_{S: j \notin S} \int_{E^S} J_S(x, dz_S)(f(x^{z_S}) - f(x)),$$

we obtain, by summarizing the previous discussions,

$$D_{x_j \to y_j} \mathcal{L} f(x) = \mathcal{L}_j [D_{x_j \to y_j} f](x)$$

(3.6)
$$- \sum_{S: j \in S} J(x, E^S) D_{x_j \to y_j} f + \sum_S g_{Sj}$$

$$|g_{Sj}| \le c_S(j) d(x_j, y_j) \sum_{i \in S} \delta_i(f).$$

Dividing both sides by $d(x_j, y_j)$ and noting

(3.7) $$\frac{\mathcal{L}_j g(x)}{d(x_j, y_j)} = \left(\mathcal{L}_j \frac{g}{d(\cdot_j, y_j)}\right)(x)$$

(with the convention that $0/0 := 0$), we get

$$\delta_j(\mathcal{L}f) \le \|\mathcal{L}_j\|_{b\mathcal{F}_T} \delta_j(f)$$
$$+ \sup_x \sum_{S: j \in S} J(x, E^S) \delta_j(f) + \sum_S c_S(j) d(x_j, y_j) \sum_{i \in S} \delta_i(f).$$

Hence, $f \to \mathcal{L}f$ is bounded on $C_{\text{Lip}}$.

*Step* 2. Note that $\mathcal{L}_j f := \sum_{S: j \notin S} \int_{E^S} J_S(x, dz_S)(f(x^{z_S}) - f(x))$ generates again a Markov semigroup on $b\mathcal{F}_T$. Let $V_j(x) := \sum_{S: j \in S} J(x, E^S)$. By the definition of $\eta$ in (3.4), $V_j \ge \eta$. Consequently, by the Feynman–Kac formula, $(\lambda + V_j - \lambda \mathcal{L}_j)^{-1}$ is bounded, positive on $b\mathcal{F}_T$ with norm bounded by $(\lambda + \eta)^{-1}$ for any $\lambda > 0$ (this is the key for the proof).

For any $\lambda > 0$ and $f \in C_{\text{Lip}}(E^T)$, let $g := f - \lambda \mathcal{L}f$. By (3.6), we have

$$(\lambda + V_j) D_{x_j \to y_j} f - \lambda \mathcal{L}_j D_{x_j \to y_j} f \le D_{x_j \to y_j} g + \lambda \sum_S g_{Sj}.$$

Dividing both sides above by $d(x_j, y_j)$ and putting $h_j(x) := \frac{D_{x_j \to y_j} f(x)}{d(x_j, y_j)}$, we obtain, by (3.7) and the estimate (3.6) for $g_{Sj}$,

$$(\lambda + V_j - \mathcal{L}_j) h_j \le \delta_j(g) + \lambda \sum_S c_S(j) \sum_{i \in S} \delta_i(f) = \delta_j(g) + \lambda \sum_i c_{ij} \delta_i(f).$$

Since the resolvent $(\lambda + V_j - \lambda \mathcal{L}_j)^{-1}$ on $b\mathcal{F}_T$ is positive with norm bounded by $(\lambda + \eta)^{-1}$ as noted above, we get

$$\sup_x h_j(x) \le \frac{1}{\lambda + \eta}\left(\delta_j(g) + \lambda \sum_i c_{ij} \delta_i(f)\right).$$

Since $y_j$ is arbitrary, we get

$$\delta_j(f) \le \frac{1}{\lambda + \eta}\left(\delta_j(g) + \lambda \sum_i c_{ij} \delta_i(f)\right).$$

*Step* 3. The rest of the proof is same as that of Proposition 2.5. □

On the product space, consider the metric

(3.8) $$d_{l_1}(x, y) := \sum_{i \in T} d(x_i, y_i).$$

It is the usual Hamming distance if $d(x,y) = \mathbb{1}_{x \neq y}$. Theorem 3.1 yields an explicit estimate of the exponential decay of $P_t$ in $C_{\text{Lip}}$. Let us translate it as (this type of translation has been given by Zhang [39]) the following:

COROLLARY 3.5. *In the context of Theorem 3.1, assume $r_{\text{sp}}(C) < 1$. Then $(P_t)$ has a unique invariant measure $\mu$ such that $\int_{E^T} d_{l_1}(x,y)\mu(dy) < +\infty$ for every (or some) $x$; moreover, for each $x \in E^T$,*

$$W_{1,d_{l_1}}(P_t(x,\cdot),\mu) \leq e^{-\eta t} \max_j \sum_i (e^{tC})_{ij} \int_{E^T} d_{l_1}(x,y)\mu(dy)$$

$$\leq e^{-t(\eta - \|C\|_1)} \int_{E^T} d_{l_1}(x,y)\mu(dy).$$

PROOF. Notice that the Lipschitzian coefficient $\|f\|_{\text{Lip}(d_{l_1})}$ of $f$ w.r.t. $d_{l^1}$ equals exactly to $\max_{i \in T} \delta_i(f)$. Hence, for any $f \in C_{\text{Lip}}(E^T)$ such that $\|f\|_{\text{Lip}(d_{l_1})} = \max_{i \in T} \delta_i(f) \leq 1$, we have, by (3.5),

$$\max_{j \in T} \delta_j(P_t f) \leq \max_{j \in T} e^{-\eta t} \sum_i (e^{tC})_{ij} = e^{-\eta t} \|e^{tC}\|_1$$

$$\leq e^{-\eta t} e^{t\|C\|_1} = e^{-t(\eta - \|C\|_1)}.$$

Thus, for every $x, y \in E^T$,

$$|P_t f(x) - P_t f(y)| \leq d_{l_1}(x,y) e^{-\eta t} \|e^{tC}\|_1.$$

Let $\nu_1, \nu_2 \in M_1^{d_{l_1},1}(E^T)$ and $\pi(dx,dy)$ a coupling of $\nu_1, \nu_2$. In the inequality above, integrating w.r.t. $\pi(dx,dy)$, and next taking the infimum over all couplings $\pi(dx,dy)$ of $(\nu_1,\nu_2)$, we get

$$|(\nu_1 P_t)f - (\nu_1 P_t)f| \leq W_{1,d_{l_1}}(\nu_1,\nu_2) e^{-\eta t} \|e^{tC}\|_1,$$

where it follows, by (2.2),

(3.9) $$W_{1,d_{l_1}}(\nu_1 P_t, \nu_2 P_t) \leq W_{1,d_{l_1}}(\nu_1,\nu_2) e^{-\eta t} \|e^{tC}\|_1.$$

As the last quantity tends to zero (exponentially), there is some $t_0 > 0$ such that $\nu \to \nu P_t$ is a strict contraction on the complete metric space $(M_1^{d_{l_1},1}(E^T), W_{1,d_{l_1}})$ [32]. Hence, $P_{t_0}$ has a unique invariant measure $\mu \in M_1^{d_{l_1},1}$ by the fixed point theorem. Moreover, the previous relation implies that $P_{nt_0}(x,\cdot) \to \mu$ in the Wasserstein distance for every $x \in E^T$. Thus, $\mu$ is the unique invariant measure of $P_{t_0}$, therefore that of $(P_t)$.

Now the desired estimate follows immediately by (3.9). $\square$

## 4. Transportation inequality $T_1$ for Gibbs measures.

4.1. *An interpretation of the a priori estimate of Dobrushin.* We begin with the fundamental a priori estimate of Dobrushin [9].

LEMMA 4.1. *Let $\mu$ be the Gibbs measure associated with the given one-point specification $(\mu_i)_{i \in T}$. Assume that the spectral radius $r_{\rm sp}(C)$ of the Dobrushin interdependence matrix is strictly smaller than* 1. *Then for any probability measure $\nu$ on $E^T$ [denoted by $\nu \in M_1(E^T)$] and $f \in C_{\rm Lip}(E^T)$,*

$$(4.1) \qquad \left| \int_{E^T} f \, d(\mu - \nu) \right| \leq \sum_{i,j} \delta_i(f) D_{ij} \mathbb{E}^\nu W_1^d(\mu_j, \nu_j),$$

*where $D := (I - C)^{-1} = \sum_{n=0}^\infty C^n$, $\nu_i = \nu(\cdot / \mathcal{F}_{T \setminus \{i\}})$ (regular conditional distribution of $x_i$ under $\nu$ knowing $\mathcal{F}_{T \setminus \{i\}}$).*

This is due to Dobrushin [9] when $d$ is the trivial metric, and is extended to general metrics by Föllmer [12] (see also [13], Theorem 8.20).

Consider the metric $d_{l_1}(x, y)$ given in (3.8). Notice that the Lipschitzian coefficient $\|f\|_{\text{Lip}(d_{l_1})}$ of $f$ w.r.t. $d_{l^1}$ equals exactly to $\max_{i \in T} \delta_i(f)$. Hence, taking the supremum in (4.1) over all $f$ such that $\max_{i \in T} \delta_i(f) \leq 1$, we obtain, by (2.2),

$$W_1^{d_{l_1}}(\mu, \nu) \leq \sum_{i,j} D_{ij} \mathbb{E}^\nu W_1^d(\mu_j, \nu_j) \leq \sup_j \sum_i D_{ij} \sum_j \mathbb{E}^\nu W_1^d(\mu_j, \nu_j).$$

Furthermore, it is obvious that if $\|C\|_1 = \sup_j \sum_i C_{ij} < 1$,

$$\sup_j \sum_i D_{ij} = \|(I - C)^{-1}\|_1 \leq \frac{1}{1 - \|C\|_1}.$$

Consequently, we have shown the following:

PROPOSITION 4.2. *Assume that $\|C\|_1 < 1$. Then for any probability measure $\nu$ on $E^T$,*

$$(4.2) \qquad W_1^{d_{l_1}}(\mu, \nu) \leq \frac{1}{1 - \|C\|_1} \sum_j \mathbb{E}^\nu W_1^d(\mu_j, \nu_j).$$

This result is the counterpart for the Wasserstein distance of Theorem 2.1.

4.2. *$T_1$-transportation inequality and Hoeffding's inequality.* Now recall that $\mu$ is said to satisfy the $T_1$-transportation inequality w.r.t. the metric $d$, if there is some constant positive $K$ such that, for all probability measures $\nu$,

$$(4.3) \qquad W_1^d(\nu, \mu) \leq \sqrt{2K h(\nu/\mu)},$$

where $h(\nu/\mu)$ is the relative entropy (or Kullback information) of $\nu$ w.r.t. $\mu$, given by

$$h(\nu/\mu) := \begin{cases} \int \frac{d\nu}{d\mu} \log \frac{d\nu}{d\mu} \, d\mu, & \text{if } \nu \ll \mu, \\ +\infty, & \text{otherwise.} \end{cases}$$

This relation will be denoted by $\mu \in T_1(K/d)$. Recall that when $d$ is the trivial metric, (4.3) holds with $K = 1/4$, which is the well-known Pinsky–Csiszär inequality.

In the following result we assume, moreover, that, for each $S \subset T$, there is a regular conditional distribution $\mu_S(dx_S|x)$ of $x_S$ knowing $x_{T\setminus S}$ under $\mu$ such that, for each $i \in S$, $\mu_i(\cdot/x)$ (fixed at the beginning) constitutes a regular conditional distribution of $x_i$ knowing $x_{S\setminus\{i\}}$ under $\mu_S(\cdot/x)$ for *every* $x \in E^T$.

THEOREM 4.3. *Assume $\|C\|_1 < 1$ and that there is some constant $K > 0$ such that $\mu$-a.s.,*

(4.4) $\qquad \mu_S(dx_i|x) \in T_1(K/d) \qquad \forall\, x \in E^T, i \in S \subset T.$

*Then for any probability measure $\nu$ on $E^T$,*

(4.5) $\qquad W_1^{d_{l_1}}(\nu, \mu) \leq \sqrt{\frac{2K|T|}{(1 - \|C\|_1)^2} h(\nu/\mu)},$

*that is, $\mu \in T_1(K(1-\|C\|_1)^{-2}|T|/d_{l_1})$. Equivalently (due to Bobkov–Götze), for any $F : E^T \to \mathbb{R}$ such that $\max_{i\in T} \delta_i(F) = \alpha < +\infty$,*

(4.6) $\qquad \mathbb{E}^\mu \exp(F - \mathbb{E}^\mu F) \leq \exp\left(\frac{K|T|\alpha^2}{2(1-\|C\|_1)^2}\right),$

*where $|T|$ is the cardinality of $T$.*

*In particular, when the diameter of $E$, $D := \sup_{x,y\in E} d(x,y) < +\infty$, both (4.5) and (4.6) hold with $K = D^2/4$.*

Before proving this theorem, let us give a quick application. Assume that $f : E \to \mathbb{R}$ is a bounded measurable function with $a \leq f \leq b$. Consider the functional related with the CLT,

$$F(x) := \frac{1}{\sqrt{|T|}} \sum_{i \in T} (f(x_i) - \mathbb{E}^\mu f(x_i)).$$

Then w.r.t. the trivial metric $d$, $\delta_i(F) \leq (b-a)/\sqrt{|T|}$ for every $i \in T$. Thus, by Theorem 4.3, (4.6),

$$\mathbb{E}^\mu e^{\lambda F} \leq \exp\left(\frac{1}{8(1-\|C\|_1)^2} \lambda^2 (b-a)^2\right) \qquad \forall\, \lambda \in \mathbb{R},$$

which, in the independent case, is the well-known sharp Hoeffding inequality.

REMARK 4.4. For an $E$-valued homogeneous Markov chain $(X_k)_{k\geq 1}$ with transition kernel $P(x,dy)$, Marton [24] proved that the law $\mathbb{P}_n$ of $(X_k)_{1\leq k\leq n}$ satisfies the transportation inequality "$T_1$" w.r.t. the Hamming metric on $E^n$ with the constant $K_n = \frac{n}{4(1-r)^2}$, where

$$r := \tfrac{1}{2}\sup_{x\neq y}\|P(x,\cdot) - P(y,\cdot)\|_{\mathrm{TV}} = W_{d,1}(P(x,\cdot), P(y,\cdot))$$

($d$ being the trivial metric). This result is generalized by Djellout, Guillin and the author [7] to general stochastic sequences w.r.t. a general metric. One can then regard Theorem 4.3 as a generalization of those results to the case of random fields.

After the first version of this paper was submitted, we learned a new work of Marton [25] in which she establishes the $T_2$-transportation inequality for Gibbs measures, by means of similar Dobrushin's interpendence coefficients (but her approach is completely different). Her $T_2$-transportation inequality, though qualitatively stronger than the $T_1$'s in Theorem 4.3, contains, however, an extra absolute constant (it is then much less precise) and works only in continuous spin space cases.

The study on transportation inequalities is very active at present; see Villani [32] and the recent thesis of Gozlan [15] for an account of art. On a Riemannian manifold, the $T_1$-transportation inequality is weaker than the log-Sobolev inequality, but is neither stronger nor weaker than the Poincaré inequality. For instance, the measure $e^{-|x|}\,dx/2$ satisfies the Poincaré inequality, but not the $T_1$-transportation inequality; the measure $(\mathbf{1}_{[-2,-1]}(x)+\mathbf{1}_{[1,2]}(x))\,dx/2$ satisfies the $T_1$-transportation inequality, but not Poincaré's.

PROOF OF THEOREM 4.3. We prove at first a general known claim: any probability measure $\nu$ on $(E,d)$ with $D = \sup_{x,y\in E} d(x,y) < +\infty$ satisfies $T_1(K/d)$ with $K = D^2/4$. In fact, for any $f \in C_{\mathrm{Lip}}(E)$ with $\|f\|_{\mathrm{Lip}} \leq 1$, $\delta(F) = \sup_x F(x) - \inf_x F(x) \leq D$. Hence, we have (a good exercise for undergraduate students)

$$\mathbb{E}^\nu e^{f-\nu(f)} \leq e^{D^2/8}.$$

This implies the desired $T_1(K/d)$ with $K = D^2/4$ by Bobkov–Götze's theorem [2].

The equivalence between (4.5) and (4.6) follows from Bobkov–Götze's theorem (cf. [2]) and the fact that the Lipschitzian coefficient $\|F\|_{\mathrm{Lip}(d_{l_1})}$ of $F$ w.r.t. $d_{l^1}$ equals exactly to $\max_{i\in T} \delta_i(F)$.

Let us prove (4.6) by the martingale method (as in [7]) in two steps.

*Step* 1. Identifying $T$ as $\{1, 2, \ldots, n\}$, where $n = |T|$ (the cardinality of $T$), we consider the martingale

$$M_0 = \mathbb{E}^\mu F, \qquad M_k(x_1^k) = \int F(x_1^k, x_{k+1}^n) \mu(dx_{k+1}^n | x_1^k), \qquad i \geq 1,$$

where $x_i^j = (x_i, x_{i+1}, \ldots, x_j)$, $\mu(dx_{k+1}^n | x_1^k) = \mu_S(dx_{k+1}^n | x)$ with $S = \{k+1, \ldots, n\}$, given previously. Since $M_n = F$, we have

$$\mathbb{E}^\mu e^{F - \mathbb{E}^\mu F} = \mathbb{E}^\mu \exp\left(\sum_{k=1}^n (M_k - M_{k-1})\right).$$

By recurrence, for (4.6), it suffices to establish that, for each $k = 1, \ldots, n$, $\mu$-a.s.,

(4.7)
$$\int \exp(M_k(x_1^{k-1}, x_k) - M_{k-1}(x_1^{k-1})) \mu(dx_k / x_1^{k-1})$$
$$\leq \exp\left(\frac{K\alpha^2}{2(1 - \|C\|_1)^2}\right).$$

*Step* 2. By the assumption (4.4) and the Bobkov–Götze theorem, for (4.7), it is enough to show that

$$|M_k(x_1^{k-1}, x_k) - M_k(x_1^{k-1}, y_k)| \leq \frac{\alpha}{1 - \|C\|_1} d(x_k, y_k).$$

By the triangle inequality, we have

$$|M_k(x_1^{k-1}, x_k) - M_k(x_1^{k-1}, y_k)|$$
$$= \left|\int F(x) \mu(dx_{k+1}^n / x_1^k) - \int F(x^{y_k}) \mu(dx_{k+1}^n / x_1^{k-1}, y_k)\right|$$
$$\leq \left|\int F(x) - F(x^{y_k}) \mu(dx_{k+1}^n / x_1^{k-1}, y_k)\right|$$
$$+ \left|\int F(x) [\mu(dx_{k+1}^n / x_1^k) - \mu(dx_{k+1}^n / x_1^{k-1}, y_k)]\right|$$
$$\leq \alpha d(x_k, y_k) + \left|\int F(x) [\mu(dx_{k+1}^n / x_1^k) - \mu(dx_{k+1}^n / x_1^{k-1}, y_k)]\right|.$$

By the dual characterization (2.2), the last term above is

$$\leq \alpha W_1^{d_{l_1}}(\mu(dx_{k+1}^n / x_1^{k-1}, y_k), \mu(dx_{k+1}^n / x_1^k)).$$

Now by Proposition 4.2, this quantity is bounded from above by

$$\frac{\alpha}{1 - \|C\|_1} \sum_{l=k+1}^n \mathbb{E}^{\mu(\cdot / x_1^{k-1}, y_k)} W_1^d(\mu_l(\cdot / x), \mu_l(\cdot / x^{y_k}))$$
$$\leq \frac{\alpha}{1 - \|C\|_1} \sum_{l=k+1}^n c_{lk} d(x_k, y_k) \leq \frac{\alpha \|C\|_1}{1 - \|C\|_1} d(x_k, y_k).$$

Thus, in summary we have

$$|M_k(x_1^{k-1}, x_k) - M_k(x_1^{k-1}, y_k)| \leq \left(\alpha + \frac{\alpha \|C\|_1}{1 - \|C\|_1}\right) d(x_k, y_k)$$
$$= \frac{\alpha}{1 - \|C\|_1} d(x_k, y_k),$$

the desired estimate. □

**5. Several concrete examples.** In this section we consider Gibbs measures on $E^{\mathbb{Z}^d}$ associated with interaction $(\phi_S)_{S \subset\subset \mathbb{Z}^d}$, where $S \subset\subset \mathbb{Z}^d$ means that $S$ is a finite subset of $\mathbb{Z}^d$. More precisely, for each finite $T \subset \mathbb{Z}^d$ and the boundary condition $y$ on $T^c$, the (local) Gibbs measure $\mu_T(dx_T|y)$ is given by

$$\frac{\exp(\sum_{S \cap T \neq \varnothing} \phi_S(x_T y_{T^c}))}{Z_T(y)} \prod_{i \in T} m(dx_i),$$

where $m$ is some reference $\sigma$-finite measure, and $Z_T(y)$ is the normalization constant. Here $\phi_S$ is $\mathcal{F}_S := \sigma(x_i, i \in S)$-measurable. Though the Dobrushin uniqueness condition works for general graphs rather than the lattice $\mathbb{Z}^d$ and all our results below have easy counterparts for graphs, we choose still the lattice $\mathbb{Z}^d$ because in that case the known results are numerous and the reader could compare more easily then with ours.

5.1. *Two examples of discrete spin models.* For a wide variety of discrete spin models w.r.t. the discrete metric $d(x,y) = \mathbb{1}_{x \neq y}$, the Dobrushin interdependence matrix $C$ has been estimated explicitly; see [13] and references therein.

EXAMPLE 5.1. $E := \{-1, 1\}, \Phi_S(x) := -J(S)x^S$, where $x^S := \prod_{i \in S} x_i$, $m$ is the counting measure on $\{-1, 1\}$. For this model, by Georgii [13], page 145, for each $\mu = \mu_T(\cdot|x)$, where $T$ is a finite subset of $\mathbb{Z}^d$,

$$\max\{\|C\|_1, \|C\|_\infty\} \leq \sup_{i \in \mathbb{Z}^d} \sum_{S : S \ni i} (|S| - 1) \tanh |J(S)| =: r.$$

(This estimate is optimal in some sense.) So when $r < 1$, all our general results apply.

When $d = 1$, $J(S) = \beta$ if $S = \{i, j\}$ with $|i - j| = 1$ and $J(S) = 0$, otherwise (one-dimensional nearest-neighbor Ising model), a pretty result of Minlos and Trishch [27] says that

$$\lambda_1(\mu) = 1 - \tanh \beta.$$

Theorem 2.1 yields only $\lambda_1(\mu) \geq 1 - r = 1 - 2 \tanh \beta$. Anyway, as there is no phase transition in the one-dimensional case, certainly one should use other parameters than the Dobrushin interdependence matrix to yield an explicit estimate of $\lambda_1(\mu)$ in dimension one.

Of course, our results become interesting when $d \geq 2$.

EXAMPLE 5.2 (*Potts anti-ferromagnet for large number of spin states*). $E = \{1, 2, \ldots, N\}$, $m$ is the counting measure on $E$ and

$$\Phi_S(x) = \begin{cases} J\mathbb{1}_{x_i = x_j}, & \text{if } S = \{i, j\}, \ |i - j| = 1, \\ 0, & \text{otherwise.} \end{cases}$$

Here $J > 0$ is a constant. For this model, Salas and Sokal [28] (communicated to the author by one referee) proved that w.r.t. the trivial metric, $c_{ij} \leq (N - 2d)^{-1}$ for $|i - j| = 1$. Hence, for every $\mu = \mu_T(dx_T|x)$ where $T$ is a finite subset of $\mathbb{Z}^d$,

$$\max\{\|C\|_1, \|C\|_\infty\} \leq \frac{2d}{N - 2d},$$

independent of the interaction strength $J$. Once $N > 4d$, all our general results apply.

For this example, Proposition 2.5 is much better than Liggett's $M - \varepsilon$ theorem.

5.2. *Two continuous spin models.*

EXAMPLE 5.3 (*N-vector model* [14, 20, 33]). Let $E = S^p$ ($p \geq 1$ integer), the unit sphere in $\mathbb{R}^{p+1}$, equipped with the normalized Lebesgue measure $m(dx)$, and

$$\Phi_S(x) = \begin{cases} -J(i - j)x_i \cdot x_j, & \text{if } S = \{i, j\} (i \neq j), \\ 0, & \text{otherwise,} \end{cases}$$

where $x \cdot y$ is the standard inner product in $\mathbb{R}^{p+1}$, and the interaction coefficients $\{J(i)\}_{i \in \mathbb{Z}^d}$ with $J(0) = 0$ is pair and absolutely summable, that is, $J(-i) = J(i)$ and

$$\gamma := \sum_{i \in \mathbb{Z}^d} |J(i)| < +\infty.$$

This is the so-called $N$-vector model with $N = p + 1$.

We begin with the estimate of $\lambda_0$ in Theorem 2.2 and of the constant $K$ in the transportation inequality of $\mu_T(dx_i|x)$ in Theorem 4.3. For every $h \in \mathbb{R}^{p+1}$, consider the probability measure on $S^p$,

(5.1) $$\mu_h(dx) := \frac{1}{Z(h)} e^{h \cdot x} m(dx),$$

where $Z(h)$ is the normalization constant.

LEMMA 5.4. *Let $\nabla$ be the Riemannian gradient on $S^p$, then*

$$(5.2) \quad \lambda_0(p,\gamma) := \inf_{h \in \mathbb{R}^{p+1}, |h| \leq \gamma} \lambda_1(\mu_h, \nabla) \geq \frac{(p-1-\gamma)\pi^2}{8(1-\exp[-(\pi^2/8)(p-1-\gamma)])}.$$

*In particular, for all $i \in \mathbb{Z}^d, x \in (S^p)^{\mathbb{Z}^d}$,*

$$(5.3) \quad \lambda_1(\mu_i(dx_i|x), \nabla) \geq \lambda_0(p,\gamma).$$

*Furthermore, w.r.t. the Riemannian metric $d$ on $S^p$, for every finite $T \subset \mathbb{Z}^d$,*

$$(5.4) \quad \mu_T(dx_i|x) \in T_1(K_0/d), \qquad K_0 = \min\left\{\frac{\pi^2}{4}; \frac{e^{2\gamma}}{p}\right\}.$$

PROOF. By the famous Lichnerowicz estimate [4], we have $\lambda_1(\mu_0, \nabla) = p$. Now for every measure $\nu(dx) = e^{-W}m(dx)/C$, where $W \in C^2(S^p)$ whose Hessian matrix $HessW \geq \beta I$ (in the order of nonnegative definiteness) and $C$ is the normalization constant, the Bakry–Emery curvature [1] $Ric(S^p) + Hess(W)$ of $\nu$ is bounded from below by $Ric(S^p) + \beta = p - 1 + \beta$, we have

$$\lambda_1(\nu, \nabla) \geq \frac{(p-1-\beta)\pi^2}{8(1-\exp[-(\pi^2/8)(p-1-\beta)])},$$

by a sharp estimate due to Chen and Wang [6]. Now for $\nu = \mu_h$, we have $W(x) = hx$ and $HessW \geq -|h|I \geq -\gamma I$ whenever $|h| \leq \gamma$, so we get (5.2). Hence, (5.3) follows since $\mu_i(dx_i|x) = \mu_h(dx_i)$ with

$$h = \sum_{j : j \neq i} J(i-j)x_j, \qquad |h| \leq \gamma.$$

For (5.4), by Theorem 4.3, we have $K_0 \leq \pi^2/4$. To show $K_0 \leq e^{2\gamma}/p$, we begin with the following sharp log-Sobolev inequality due to Bakry–Emery [1]:

$$Ent_m(f^2) \leq \frac{2}{p}\int_{S^p} |\nabla f|^2 m(dx) \qquad \forall f \in C^1(S^p),$$

where $Ent_\nu(f) = \mathbb{E}^\nu f \log \frac{f}{\nu(f)}$ is the relative entropy of $f \geq 0$ w.r.t. $\nu$. Since the marginal law of $\mu_T(\cdot|x)$ at $x_i$,

$$\nu_{T,i}(dx_i) := \mu_T(dx_i|x) = \int \mu_i(dx_i|x)\mu_T(dx_{T\setminus\{i\}}|x) = e^{-W(x_i)}m(dx_i)/C,$$

where $W$ satisfies $\delta(W) := \sup_{x_i} W - \inf_{x_i} W \leq 2\gamma$, hence, $\nu_{T,i}$ satisfies the log-Sobolev inequality,

$$Ent_{\nu_{T,i}}(f^2) \leq 2\frac{e^{2\gamma}}{p}\int_{S^p} |\nabla f|^2 \, d\nu_{T,i} \qquad \forall f \in C^1(S^p),$$

by a remark in [19]. This implies, by Ledoux [18] (using the Herbst method),

$$\mathbb{E}^{\nu_{T,i}} e^{f - \nu_{T,i}(f)} \leq \exp\left(\frac{e^{2\gamma} \|f\|_{\text{Lip}}^2}{2p}\right) \quad \forall f \in C^1(S^p),$$

which is equivalent to (Bobkov–Götze's theorem)

$$\nu_{T,i} \in T_1(K/d), \qquad K = \frac{e^{2\gamma}}{p}.$$

Thus, (5.4) is established. □

PROPOSITION 5.5. *W.r.t. the Euclidean metric $d_E$ of $\mathbb{R}^{p+1}$ restricted to $S^p$, the coefficient of interdependence of Dobrushin satisfies*

$$(5.5) \qquad c_{ij}^E \leq |J(i-j)| \frac{\sigma_E(p,\gamma)}{\sqrt{p+1}},$$

*where*

$$(5.6) \qquad \sigma_E^2(p,\gamma) := \sup_{f,h} \mu_h(f;f) \leq \min\{1, 1/\lambda_0(p,\gamma)\}.$$

*Here the supremum is taken over all $h \in \mathbb{R}^{p+1}$ with $|h| \leq \gamma$ and $f : S^p \to \mathbb{R}$ such that its Lipschitzian coefficient $\|f\|_{\text{Lip}(d_E)}$ w.r.t. $d_E$ is less than 1, and $\lambda_0(p,\gamma)$ is given in (5.3). In particular, if*

$$(5.7) \qquad \gamma := \sum_{i \in \mathbb{Z}^d} |J(i)| < \frac{\sqrt{p+1}}{\sigma_E(p,\gamma)},$$

*then for every finite $T \subset \mathbb{Z}^d$,*

$$(5.8) \quad \lambda_1(\mu_T) \geq 1 - \frac{\gamma \sigma_E(p,\gamma)}{\sqrt{p+1}}, \qquad \lambda_1(\mu_T, \nabla) \geq \left(1 - \frac{\gamma \sigma_E(p,\gamma)}{\sqrt{p+1}}\right) \lambda_0(p,\gamma).$$

*Furthermore, if*

$$(5.9) \qquad \gamma < \sqrt{(p+1)\lambda_0(p,\gamma)}$$

*[stronger than (5.7)], then*

$$(5.10) \quad \mu_T \in T_1(\tilde{K}|T|/d_{l_1}), \qquad \tilde{K} = \frac{\min\{e^{2\gamma}/p, \pi^2/4\} \cdot \lambda_0(p,\gamma)(p+1)}{(\sqrt{(p+1)\lambda_0(p,\gamma)} - \gamma)^2},$$

*where $d_{l_1}(x_T, y_T) := \sum_{i \in T} d(x_i, y_i)$ (d being the Riemannian metric on $S^p$).*

PROOF. For the spectral gap estimate (5.8), it is enough to prove (5.5) and (5.6) by Theorem 2.2.

Let $i, j \in \mathbb{Z}^d$ two different sites. Given $x, y \in (S^p)^{\mathbb{Z}^d}$ such that $x = y$ off $j$, let

$$h(t) := \sum_{k \neq i} J(i - k)((1-t)x_k + ty_k), \qquad t \in [0, 1],$$

and consider $\mu_{h(t)}$ as given by (5.1). Then $\mu_{h(0)} = \mu_i(\cdot | x)$ and $\mu_{h(1)} = \mu_i(\cdot | y)$. For every $f : S^p \to \mathbb{R}$ such that $\|f\|_{\text{Lip}(d_E)} \leq 1$, we have

$$
\begin{aligned}
\int_{S^p} f \, d(\mu_1 - \mu_0) &= \int_0^1 \frac{d}{dt} \int_{S^p} f \, d\mu_{h(t)} \, dt \\
&= |y_j - x_j| \int_0^1 \mu_{h(t)}(f; e \cdot x) \, dt,
\end{aligned}
\tag{5.11}
$$

where $e = (y_j - x_j)/|y_j - x_j|$. Since $|h(t)| \leq \gamma$, we have by Cauchy–Schwarz,

$$\mu_{h(t)}(f; e \cdot x) \leq \sqrt{\mu_h(f; f) \cdot \mu_h(t)(e \cdot x; e \cdot x)} \leq \sigma_E(p, \gamma) \sqrt{\mu_{h(t)}(e \cdot x; e \cdot x)}.$$

Hence, for (5.5), by the dual formula (2.2), it remains to prove that, for any $h, e \in \mathbb{R}^{p+1}$ with $|e| = 1$,

$$\mu_h(e \cdot x; e \cdot x) \leq \frac{1}{p+1}. \tag{5.12}$$

Let $\hat{e}_k, k = 1, \ldots, p+1$ be an orthonormal basis of $\mathbb{R}^{p+1}$ such that $\hat{e}_1 \cdot h = |h|$, and $e$ is a linear combination of $e_1, e_2$. Let $\hat{x}_k := \hat{e}_k \cdot x$. We have

$$\mu_h(\hat{x}_k; \hat{x}_k) = \mu(\hat{x}_k)^2 \leq \frac{1}{p+1} \qquad \forall k \geq 2,$$

by [20], (3.15) (due to Dyson–Lieb–Simon [11]) and

$$\mu_h(\hat{x}_1; \hat{x}_1) \leq \frac{1}{p+1}$$

by Levin [20], (3.28) (those two estimates can be easily derived by the correlation inequalities on $S^p$ in [14]). Now noting that $\hat{x}_1, \hat{x}_1$ are orthogonal in $L^2(S^p, \mu_h)$, we obtain (5.12) and so the desired (5.5).

To prove (5.6), note that, for any $f \in C^1(S^p)$ such that $\|f\|_{\text{Lip}(d_E)} \leq 1$, we have $a \leq f \leq b$ with two constants satisfying $b - a \leq 2$ (for the diameter of $S^p$ w.r.t. $d_E$ is 2). Thus,

$$\mu_h(f; f) \leq \int_{S^p} \left( f - \frac{a+b}{2} \right)^2 d\mu_h \leq \frac{(b-a)^2}{4} \leq 1,$$

where it follows Simon's bound $\sigma_E^2(p, \gamma) \leq 1$ [29]. Furthermore, as the Riemannian distance $d$ on $S^p$ is larger than $d_E$, we also have $\|f\|_{\text{Lip}(d)} \leq 1$. Thus, by the Poincaré inequality,

$$\mu_h(f; f) \leq \frac{1}{\lambda_1(\mu_h, \nabla)} \int_{S^p} |\nabla f|^2 \, d\mu_h \leq \frac{1}{\lambda_1(\mu_h, \nabla)},$$

which yields $\sigma_E^2(p,\gamma) \leq 1/\lambda_0(p,\gamma)$. This completes the proof of (5.6).

To prove the $T_1$-transportation inequality, we should estimate the coefficients $c_{ij}$ of interdependence of Dobrushin w.r.t. the Riemannian metric $d$. By the same proof as that of (5.5), we have

$$c_{ij} \leq |J(i-j)|\frac{\sigma_R(p,\gamma)}{\sqrt{p+1}},$$

where

$$\sigma_R^2(p,\gamma) := \sup_{\|f\|_{\mathrm{Lip}(d)}\leq 1, |h|\leq \gamma} \mu_h(f;f).$$

By the Poincaré inequality as above, we have

$$\sigma_R^2(p,\gamma) \leq \frac{1}{\lambda_0(p,\gamma)}.$$

Hence, (5.10) follows by Theorem 4.3. □

REMARK 5.6. Let us compare Proposition 5.5 with known results:

(i) For the uniqueness of Gibbs measure, "$\gamma < 1$" is Faris' condition, the better condition "$\gamma < \sqrt{p+1}$" is due to Simon [29]. And when $p \geq 5$, Levin [20] improved the bound of Simon as follows:

$$\gamma < \frac{p+1}{\sqrt{5}}.$$

Even our stronger condition (5.9) is better than Simon's for $p \geq 4$ [since $\lambda_0(p,\gamma) > 1$ once $\gamma < p-1$ by Lemma 5.4] and better than Levin's for all $p \geq 5$ since our condition (5.7) is satisfied once

(5.13) $$\gamma < \frac{2}{1+\sqrt{1+4a}}(p-1), \qquad a := \frac{8}{\pi^2} \cdot \frac{p-1}{p+1}.$$

Anyway, we are still far from the conjecture in [20], Remark 3.7, which says that the uniqueness of Gibbs' measure holds once $\gamma < p+1$, where $p+1$ is the known critical value for the phase transition of the corresponding mean field model.

(ii) Wick [33] proved that for $p = 1, 2$ and the nearest-neighbor case, if $\sqrt{\frac{5}{2}}4\gamma < 1$, then no phase transition occurs and the Glauber dynamics associated with the Dirichlet form $\mathcal{E}^\nabla$ is exponentially ergodic. His range of $\gamma$ is much more restrictive than ours.

(iii) The most simple way to obtain some explicit estimate on the spectral gap or the constant in the log-Sobolev inequality for this model is via Bakry–Emery's criterion. Indeed, since $\mu_T(dx_T|x) = e^{-H_T(x)}m^{\otimes T}(dx_T)/C_T$, where

$HessH_T \geq -2\gamma I$ on the product space $(S^p)^T$, the Bakry–Emery curvature of $\mu_T$ satisfies

$$Ric((S^p)^T) + HessW_T \geq (p-1-2\gamma)I.$$

Hence, when $2\gamma < p-1$, we have, by the criterion of Bakry–Emery [1],

$$(5.14) \quad (p-1-2\gamma)Ent_{\mu_T}(f^2) \leq 2\mathbb{E}^{\mu_T} \sum_{i \in T} |\nabla_i f|^2 \quad \forall f \in C^1((S^p)^T),$$

which implies $\lambda_1(\mu_T, \nabla) \geq p-1-2\gamma$ and $\mu_T \in T_1(K|T|/d_{l_1})$, where $K = (p-1-2\gamma)^{-1}$. Our condition (5.7) in Proposition 5.5 is better as seen from (5.13).

(iv) For this model, Zegarlinski [38], Lemma 3.2, Example 3.3, found that his coefficents of interdepence $c_{ij}^Z$ given in (2.8) verifies

$$c_{ij}^Z \leq \sqrt{\frac{2}{\lambda_0}} |J(i-j)|$$

[in his expression (3.13), $c_0$ can be replaced by $1/\lambda_0$ by the proof of Lemma 3.2 where only the single site Poincaré inequality (3.15) is used]. This estimate has an extra factor $\sqrt{2(p+1)}$ w.r.t. our estimate of $c_{ij}$ (w.r.t. the Riemannian metric) in the proof of Proposition 5.5.

Guionnet and Zegarlinski [16], basing on the Föllmer covariance estimate [12], proved the log-Sobolev inequality under the Dobrushin uniqueness condition for compact spin models, which is much stronger than $\lambda_1(\mu, \nabla) > 0$, but without a robust estimate of the involved constant.

EXAMPLE 5.7 (*The $\phi^4$ Euclidean quantum field on the lattice*). This model is given by

$$E = \mathbb{R}, \quad m(dx) = e^{-u(x)} dx/C,$$

$$\Phi_S(x) = -J(i-j)x_i x_j, \quad \text{if } S = \{i,j\} \quad \text{and} \quad \Phi_S = 0 \quad \text{otherwise},$$

where $u(x) = ax^4 - bx^2$ with $a > 0, b \in \mathbb{R}$, $C$ is the normalization constant and $J: \mathbb{Z}^d \to \mathbb{R}$ is pair and absolutely summable with $J(0) = 0$. Notice that, for every finite $T \subset \mathbb{R}^d$, every boundary condition $x \in \mathbb{R}^{\mathbb{Z}^d}$ such that $\sum_k |J(i-k)||x_k| < +\infty$ for every $i$, $\mu_T(\cdot|x)$ is well defined. In the following $\mu_T$ denotes the local Gibbs measure with such a boundary condition.

For this unbounded spin model, we can not use the trivial metric, for which $c_{ij} = +\infty$ in general. So only the Euclidean metric on $\mathbb{R}$ will be used below. We first recall a result of Helffer [17] and Ledoux [19]:

$$(5.15) \quad \lambda_1(\mu_T, \nabla) \geq \lambda_0 + h,$$

where $h$ is the infimum of the spectrum in $l^2(\mathbb{Z}^d)$ of the matrix $(\gamma_{ij})_{i,j\in\mathbb{Z}^d}$, where $\gamma_{ij} = -J(i-j)$ and

(5.16) $$\lambda_0 = \inf_{\theta \in \mathbb{R}} \lambda_1(m_\theta, \nabla), \qquad m_\theta(dx) := e^{-u(x)+\theta x} \, dx / C_\theta.$$

Applying the previous general results, we will get the following:

PROPOSITION 5.8. *Let*

(5.17) $$\sigma^2 = \text{the variance of } x \text{ under } m = \int_{\mathbb{R}} x^2 \, dm(x).$$

*If*

(5.18) $$\gamma := \sum_{k \in \mathbb{Z}^d} |J(k)| < \frac{1}{\sigma^2},$$

*then for every finite $T \subset \mathbb{Z}^d$,*

(5.19) $$\lambda_1(\mu_T) \geq 1 - \sigma^2 \gamma, \qquad \lambda_1(\mu_T, \nabla) \geq (1 - \sigma^2 \gamma) \lambda_0,$$

*where $\lambda_0$ is given in* (5.16), *and*

(5.20) $$\mu_T \in T_1(\tilde{K}|T|/d_{l_1}), \qquad \tilde{K} = \frac{K_0}{(1 - \gamma \sigma^2)^2},$$

*where $K_0$ is the best positive constant such that*

$$\mu_S(dx_i|x) \in T_1(K_0/d) \qquad \forall\, i \in S \subset\subset \mathbb{Z}^d, x.$$

REMARK 5.9. By Cassandro, Olivieri, Pellegrinotti and Presutti [3], $1/\sigma^2$ is the critical value of the interaction strength $\gamma$ of the corresponding mean field ferromagnetic model. In other words, our condition (5.18) is sharp in this point of view.

By the definition of the spectral gap,

$$\lambda_1(m, \nabla) m(x; x) = \lambda_1(m, \nabla) \sigma^2 \leq 1,$$

and since $f(x) = x$ is not an eigenfunction of the generator $\mathcal{L} = \frac{d^2}{dx^2} - u'(x) \, dx$ associated to the Dirichlet form $\int_{\mathbb{R}} (f')^2 m(dx)$, we have

$$\lambda_1(m, \nabla) \sigma^2 < 1.$$

In particular, $\lambda_0 \sigma^2 < 1$, where $\lambda_0$ is given in (5.16).

On the other hand, in the ferromagnetic case, that is, $J(i) \geq 0$, it is easy to see that the infimum $h$ of the spectrum in $l^2(\mathbb{Z}^d)$ of the matrix $(-J(i-j))_{i,j\in\mathbb{Z}^d}$ coincides with $-\sum_{k \in \mathbb{Z}^d} J(k) = -\gamma$. In such a situation, our estimate (5.19) is better than the known (5.15).

REMARK 5.10. Let $c_0$ be the best constant such that $\nu_{T,i} = \mu_T(dx_i|x)$ satisfies the log-Sobolev inequality

$$Ent_{\nu_{T,i}}(f^2) \leq 2c_0 \int_{\mathbb{R}} (f')^2 \, d\nu_{T,i} \qquad \forall f \in C_0^1(\mathbb{R}),$$

for all $i \in T \subset\subset \mathbb{Z}^d$ and all boundary conditions $x$. Ledoux [19] proved that $c_0 < +\infty$ (and an estimate of $c_0$). Hence, the "$T_1$" transportation constant $K_0$ in (5.20) satisfies $K_0 \leq c_0$ and the spectral gap $\lambda_0 \geq 1/c_0$.

When $J(\cdot)$ is of finite range, since Yoshida [37] has proven the equivalence between the Poincaré and log-Sobolev inequality, we have also the log-Sobolev inequality for $\mu_T$, uniformly over $T$ and the boundary condition, once if $\gamma < 1/\sigma^2$. A challenging open question is to give a robust estimate of the constant in that log-Sobolev inequality, better than Ledoux's [19].

REMARK 5.11. For this model, the results of Zegarlinski [38], Propositions 3.4, 3.6, do not apply.

PROOF OF PROPOSITION 5.8. Fix the finite subset $T$ of $\mathbb{Z}^d$ and the boundary condition $x$ such that $\sum_k |J(\cdot - k)||x_k| < +\infty$. Let us estimate the Dobrushin coefficient $c_{ij}$ associated with $\mu_T(\cdot|x)$ and the Euclidean metric $d$, where $i, j$ are two different sites in $T$. To this end, consider $y \in \mathbb{R}^{\mathbb{Z}^d}$ such that $x = y$ off $j$ and

$$\mu_0(dx_i) = \mu_i(dx_i|x) = e^{\psi_0} m(dx_0)/C_0,$$
$$\mu_1(dx_i) = \mu_i(dx_i|y) = e^{\psi_1} m(dx_0)/C_1,$$

where $\psi_0(x_i) = -\sum_k J(i-k) x_i x_k$, $\psi_1(x_i) = -\sum_k J(i-k) x_i y_k$. Let $\psi_t = \psi_0 + t(\psi_1 - \psi_0)$ and

$$\mu_t := \frac{e^{\psi_t} m(dx_i)}{\int_E e^{\psi_t} m(dx_i)}.$$

For any function $f \in C^1(\mathbb{R})$ with $\|f\|_{\text{Lip}} \leq 1$, we have as in the proof of Proposition 5.5,

$$\int_{\mathbb{R}} f \, d(\mu_1 - \mu_0) = \int_0^1 \mu_t(f; \psi_1 - \psi_0) \, dt.$$

Since $\|f\|_{\text{Lip}} \leq 1, \|\psi_1 - \psi_0\|_{\text{Lip}} = |J(i-j)||x_j - y_j|$, we have

$$\mu_t(f; \psi_1 - \psi_0) = \frac{1}{2} \iint_{\mathbb{R}^2} (f(x_i) - f(z_i))$$
$$\times [(\psi_1 - \psi_0)(x_i) - (\psi_1 - \psi_0)(z_i)] \, d\mu_t(x_i) \, d\mu_t(z_i)$$
$$\leq \frac{|J(i-j)||x_j - y_j|}{2} \iint_{\mathbb{R}^2} (x_i - z_i)^2 \, d\mu_t(x_i) \, d\mu_t(z_i)$$
$$= |J(i-j)||x_j - y_j| \mu_t(x; x).$$

But $\mu_t = m_\theta$ given in (5.16) for some $\theta \in \mathbb{R}$. By the famous GHS correlation inequality ([14], Corollary 4.3.4 where the condition $\theta \geq 0$ can be removed by the symmetry $x \to -x$, because we are faced only to one site),

$$\mu_t(x;x) = m_\theta(x;x) \leq m(x;x) = \sigma^2.$$

Thus,

$$W_{1,d}(\mu_i(\cdot|x), \mu_1(\cdot|y)) \leq \sigma^2 |J(i-j)| d(x_j, y_j),$$

that is,

$$(5.21) \qquad c_{ij} \leq \sigma^2 |J(i-j)|.$$

Then $\|C\|_1 \vee \|C\|_\infty \leq \sigma^2 \gamma$. Hence, $\lambda_1(\mu_T) \geq 1 - \sigma^2 \gamma$ by Theorem 2.1.

Next, since $\mu_i(\cdot|x) = m_\theta$ for some $\theta$, we have

$$\lambda_1(\mu_i(\cdot|x), \nabla) \geq \inf_{\theta \in \mathbb{R}} \lambda_1(m_\theta, \nabla) = \lambda_0,$$

hence the second estimate in (5.19) follows by Theorem 2.2.

The last transportation inequality follows by Theorem 4.3 by the estimate of $c_{ij}$ above. $\square$

**Acknowledgments.** This work has been reported in the Workshop of Beijing Normal University, May 2004 and that of Fudan University on September 2004. I am grateful to Professors M. F. Chen and J. G. Ying for their kind invitations and warm hospitality, and especially to the first for the communication of [5] and [27].

One referee has given a long constructive and conscientious report with a review of known results on the examples in Section 5 (of which the author unaware), leading to considerable improvements in Examples 5.2 and 5.3 and in the presentation of the paper in general.

LABORATOIRE DE MATHÉMATIQUES
CNRS-UMR 6620
UNIVERSITÉ BLAISE PASCAL
63177 AUBIÈRE
FRANCE
AND
DEPARTMENT OF MATHEMATICS
WUHAN UNIVERSITY
430072 CHINA
E-MAIL: Li-Ming.Wu@math.univ-bpclermont.fr